\documentclass[11pt]{article}
\usepackage[colorlinks,
            linkcolor=blue,
            anchorcolor=green,
            citecolor=red
            ]{hyperref}
\usepackage{mathrsfs,amssymb,amsfonts}
\usepackage{cite}
\usepackage{amsmath,amssymb,latexsym,color}
\usepackage{graphicx}
\usepackage[mathscr]{eucal}
\usepackage{xypic}
 \makeatletter
    
    \newcommand{\Rmnum}[1]
    {\expandafter\@slowromancap\romannumeral #1@}
    \makeatother
\newtheorem{thm}{Theorem}[section]
\newtheorem{defin}[thm]{Definition}
\newtheorem{prop}[thm]{Proposition}
\newtheorem{lemma}[thm]{Lemma}
\newtheorem{cor}[thm]{Corollary}

\newcounter{condition}[section]
\newcounter{Algorithmnumber}

\newcounter{rownumber}

\newtheorem{constr}[thm]{Construction}
\newtheorem{example}[thm]{Example}
\newtheorem{obser}[thm]{Observation}

\newtheorem{remark}{Remark}

\newenvironment{eg}{\begin{example} \rm}{\end{example}}
\newenvironment{ob}{\begin{obser} \rm}{\end{obser}}
\newenvironment{cons}{\begin{constr} \rm}{\end{constr}}

\newcommand{\proof}{{\it Proof.\quad}}
\newcommand{\qed}{\hfill\Box\medskip}

\usepackage{CJK}
\begin{document}
\begin{CJK*}{GBK}{song}
\renewcommand{\abovewithdelims}[2]{
\genfrac{[}{]}{0pt}{}{#1}{#2}}

\title{\bf On weak metric dimension of  digraphs}

\author{ Min Feng$^{1}$\quad Kaishun Wang$^{2}$\quad Yuefeng Yang$^3$\footnote{Corresponding author.}\\
{\footnotesize   \em  $^1$School of Science, Nanjing University of Science and Technology, Nanjing, 210094, China}\\
{\footnotesize   \em $^2$Sch. Math. Sci. {\rm \&} Lab. Math. Com. Sys., Beijing Normal University, Beijing, 100875,  China}\\
{\footnotesize   \em $^3$School of Science, China University of Geosciences, Beijing 100083, China}
}
 \date{}
 \maketitle

\begin{abstract}
  Using the two way distance, we introduce  the concepts of weak metric dimension of a strongly connected digraph $\Gamma$.  We first establish lower and upper bounds for the number of arcs in $\Gamma$ by using the diameter and weak metric dimension of $\Gamma$, and characterize all digraphs attaining the lower or upper bound. Then we  study a digraph with weak metric dimension $1$ and classify all vertex-transitive digraphs having weak metric dimension $1$. Finally, all digraphs of order $n$ with weak metric dimension $n-1$ or $n-2$ are determined.


%\footnote{E-mail addresses:  fengmin@njust.edu.cn (M. Feng), wangks@bnu.edu.cn (K. Wang).}
\end{abstract}

\footnote{E-mail address: fgmn\_1998@163.com (M. Feng), wangks@bnu.edu.cn (K. Wang),\\ yangyf@cugb.edu.cn (Y. Yang).}

\section{Introduction}

A directed graph (digraph) is {\em strongly connected} if for any two vertices $x$ and $y$, there is a directed path from $x$ to $y$ in this digraph.
Throughout of this paper, we always use $\Gamma$ to denote a  finite, simple and strongly connected digraph. Let $V(\Gamma)$ and $A(\Gamma)$ be the vertex  and arc set of $\Gamma$, respectively. For $x,y\in V(\Gamma)$, the {\em distance} from $x$ to $y$ in $\Gamma$, denoted by $\partial(x,y)$, is the length of a shortest path from $x$ to $y$. The {\em two way distance} from $x$ to $y$ is defined as
$$
\tilde{\partial}(x,y)=(\partial(x,y),\partial(y,x)).
$$
A subset $\{w_{1},\ldots,w_{m}\}$ of $V(\Gamma)$ is a \emph{weakly resolving set} of $\Gamma$ if
$$
(\tilde{\partial}(w_{1},u),\ldots,\tilde{\partial}(w_{m},u))\neq     (\tilde{\partial}(w_{1},v),\ldots,\tilde{\partial}(w_{m},v))
$$
for any distinct vertices $u,v\in V(\Gamma)$. The
\emph{weak metric dimension} of $\Gamma$, denoted by $\dim(\Gamma)$, is the minimum
cardinality of a weakly resolving set of $\Gamma$.

In the 1970s, the metric dimension of a graph was first introduced, by Harary and Melter \cite{HM76} and, independently, by Slater \cite{S75}. Slater referred to metric dimension of a graph as its location number. Subsequently, this parameter has appeared in various applications, as diverse as network discovery and verification \cite{Bee}, strategies for the Mastermind game \cite{Chv},  combinatorial optimization \cite{Se}
and so on. It was noted in  \cite{Kh} that determining the metric dimension of a graph is an NP-complete problem.
For more information, we refer to \cite{BC11} and \cite{CH07}.

As an analogue of the metric dimension of a graph, the directed distance dimension of an oriented graph is defined by Chartrand, Rains and  Zhang  \cite{Char1,Char2}. They defined a resolving set of a digraph $\Gamma$ as a subset $\{w_{1},\ldots,w_{m}\}$ of $V(\Gamma)$ such that
$$
(\partial(u,w_{1}),\ldots,\partial(u,w_{m}))\neq
(\partial(v,w_{1}),\ldots,\partial(v,w_{m}))
$$
for any distinct vertices $u,v\in V(\Gamma)$, and then defined the
directed distance dimension of $\Gamma$ as the minimum cardinality of a resolving set of $\Gamma$.
The directed distance dimension of a digraph $\Gamma$ was called the metric dimension of $\Gamma$ by Fehr et al. \cite{Fe}, Feng et al. \cite{fengxuwang} and  Vetr\'ik \cite{vetrik}. Bountin  et al.  \cite{boutingp}, and Abs and Vetr\'ik \cite{abasv} defined the metric dimension of a digraph $\Gamma$ as the directed distance dimension of $\Gamma^{-1}$, where $\Gamma^{-1}$ is the digraph obtained from $\Gamma$ by reversing directions of all arcs.

Relative positions of  two nodes in a directed network are determined by their two way distance. As resolving sets are often used to locate nodes, it is meaningful to use two way distance to define another type of resolving sets. Motivated by this, we have the definitions of weakly resolving sets and weak metric dimension of a digraph at the beginning of this paper.

Let $\Gamma$ be a digraph. For a subset $S$ of $V(\Gamma)$, we say that $S$ {\em weakly resolves}, or {\em resolves} for simplify, vertices $x$ and $y$, or $\{x,y\}$, if there is a vertex $z\in S$ such that $\tilde\partial(z,x)\neq\tilde\partial(z,y)$. If $\{z\}$ resolves $x$ and $y$, we also say that $z$ resolves $x$ and $y$.
A subset $S$ of $V(\Gamma)$ is a weakly resolving set of $\Gamma$ if and only if it resolves each pair of vertices in $V(\Gamma)\setminus S$.  We say that $\Gamma$ is {\em $k$-dimensional} if $\dim(\Gamma)=k$. A {\em weak metric basis}, or {\em basis}, of a $k$-dimensional digraph is a weakly resolving set of cardinality $k$.

The paper is organized as follows. In Section~2, we first use the diameter and the order of a digraph $\Gamma$ to give sharp lower and upper bounds for $\dim(\Gamma)$, and then establish lower and upper bounds for $|A(\Gamma)|$ by using the diameter and the weak metric dimension of $\Gamma$. As a result, we characterize all graphs $\Gamma$ such that $|A(\Gamma)|$ attains the lower or upper bound.
In Section~3, we study $1$-dimensional digraphs and classify all $1$-dimensional vertex-transitive digraphs. Section~4 characterizes all  $(n-2)$-dimensional digraphs of order $n$.

\section{With the diameter}

  Let $\Gamma$ be a digraph of order $n\geq 2$, then $V(\Gamma)\setminus\{v\}$ is a weakly resolving set of $\Gamma$ for each $v\in V(\Gamma)$, which implies that
  \begin{equation*}
    1\leq\dim(\Gamma)\leq n-1.
  \end{equation*}
  Actually, if we know the diameter of $\Gamma$, then we can obtain an improved upper  bound in general for $\dim(\Gamma)$,  as well as a lower bound.

\begin{thm}\label{nd}
  If $\Gamma$ is a digraph of order $n$ with diameter $d$, then
$$
f(n,d)\leq\dim(\Gamma)\leq n-d,
$$
where $f(n,d)$ is the least positive integer $k$ for which $k+d^{2k}\geq n$.
\end{thm}
\proof First, we establish the lower bound. Let $B=\{w_1,\ldots,w_k\}$ be a  basis of $\Gamma$. Then for any vertex $v\in V(\Gamma)\setminus B$, both $\partial(w_i, v)$ and $\partial(v, w_i)$ are at least $1$ and at most $d$ for each $i\in\{1,\ldots,m\}$, and so
$$
(\tilde{\partial}(w_{1},v),\ldots,\tilde{\partial}(w_{k},v))
$$
has at most $d^{2k}$  possibilities. Since $B$ is a weakly resolving set of $\Gamma$, we have $d^{2k}\geq n-k$. Thus, one gets $f(n,d)\leq k=\dim(\Gamma)$.

Next, we prove the upper bound. Let $(u_0,u_1,\ldots,u_d)$ be a shortest path from $u_0$ to $u_d$ in $\Gamma$. Write
$$
W=V(\Gamma)\setminus\{u_1,\ldots,u_d\}.
$$
Since $u_0\in W$ and $\partial(u_0,u_i)=i$ for each $i\in\{1,\ldots,d\}$, it follows that $W$ is a weakly resolving set of $\Gamma$. Hence, we have $\dim(\Gamma)\leq |W|=n-d$.
$\qed$

Note that a digraph $\Gamma$ is complete if and only if the diameter of $\Gamma$ is $1$. The following result follows immediately from Theorem~\ref{nd}.

\begin{cor}\label{complete}
  Let $\Gamma$ be a digraph of order $n\geq 2$. Then $\dim(\Gamma)=n-1$ if and only if $\Gamma$ is complete.
\end{cor}

The upper bound in Theorem~\ref{nd} is sharp. Actually, for any positive integers $n$ and $d$ with $3\leq d+1\leq n$,   we shall construct an $(n-d)$-dimensional digraph of order $n$ with diameter $d$.

\begin{eg}\label{eg1}
Given any positive integers $n$ and $d$ with $3\leq d+1\leq n$, define a digraph $\Gamma$ as follows:
$$
\begin{array}{rcl}
V(\Gamma)&=&\{v_1,v_2,\ldots,v_n\}, \\
A(\Gamma)&=&\bigcup\limits_{i=1}^{n-d+1}\{(v_n,v_i),(v_i,v_{n-d+2})\} \cup\{(v_j,v_{j+1})\mid j=n-d+2,\ldots,n-1\}.
\end{array}
$$
Then $\Gamma$ has $n$ vertices and diameter $d$. Now we show that $\dim(\Gamma)=n-d$. In fact, for any distinct vertices $v_k,v_l\in\{v_1,\ldots,v_{n-d+1}\}$ and $v_m\in V(\Gamma)\setminus\{v_k,v_l\}$, we have
$$
\tilde\partial(v_k,v_m)=\tilde\partial(v_l,v_m)=
\left\{
\begin{array}{ll}
  (d,d), &  m\leq n-d+1, \\
  (m-n+d-1,n-m+1),& m>n-d+1.
\end{array}\right.
$$
Hence, any vertex except $v_l$ and $v_k$ can not resolve $v_l$ and $v_k$. Therefore, for any basis $B$ of $\Gamma$, one gets $|B\cap\{v_k,v_l\}|\geq 1$, which implies that
$$
\dim(\Gamma)=|B|\geq|B\cap\{v_1,\ldots,v_{n-d+1}\}|\geq n-d.
$$
It follows from Theorem~\ref{nd} that $\dim(\Gamma)=n-d$.
\end{eg}

To illustrate the lower bound in Theorem~\ref{nd}, we construct a family of digraphs.

\begin{cons}\label{construction}
Given any positive integers $n$ and $d$ with $3\leq d+1\leq n$, write $k=f(n,d)$ and $m=n-f(n,d)$. Note that $m\leq d^{2k}$. For any $i\in\{1,\ldots,m\}$, there exists a unique sequence $(x_1,x_2,\ldots,x_{2k})$ with each $x_r\in\{1,\ldots,d\}$ such that
$$
i=x_1+(x_2-1)\cdot d^1+(x_3-1)\cdot d^2+\cdots+(x_{2k}-1)\cdot d^{2k-1},
$$
and then we denote
$$
v_i=(x_1,x_2,\ldots,x_{2k}).
$$
Now we define a digraph $\Gamma(n,d)$ with the vertex set
  \begin{equation}\label{uv}
  V(\Gamma(n,d))=\{u_1,\ldots,u_k\}\cup\{v_1,v_2,\ldots,v_m\},
  \end{equation}
  and there are exactly three families (F1), (F2), (F3) of arcs as follows:

    (F1) $(u_i,v_j)\in A(\Gamma(n,d))$ if and only if $a_{2i-1}=1$, where $v_j=(a_1,a_2,\ldots,a_{2k})$.

    (F2) $(v_j, u_i)\in A(\Gamma(n,d))$ if and only if $a_{2i}=1$, where $v_j=(a_1,a_2,\ldots,a_{2k})$.

    (F3) $(v_j, v_l)\in A(\Gamma(n,d))$ if and only if $v_j\ne v_l$ and
    $$
    a_{2r-1}-b_{2r-1}\geq -1\qquad\text{and}\qquad a_{2r}-b_{2r}\leq 1
    $$
    for each $r=1,\ldots,k$, where $v_j=(a_1,a_2,\ldots,a_{2k})$ and $v_l=(b_1,b_2,\ldots,b_{2k})$.
\end{cons}

Note that $|V(\Gamma(n,d))|=n$. The following result gives three families of digraphs such that the lower bound in Theorem~\ref{nd} is attained.

\begin{lemma}\label{eg2}
   {\rm(i)} For $d+1\leq n\leq d^2+1$, the diameter of  $\Gamma(n,d)$ is $d$ and
    $$
    \dim(\Gamma(n,d))=1=f(n,d).
    $$

  {\rm(ii)} For $3\leq n$, the diameter of $\Gamma(n,2)$ is $2$ and $\dim(\Gamma(n,2))=f(n,2)$.

  {\rm(iii)} For $4\leq n$, the diameter of $\Gamma(n,3)$ is $3$ and $\dim(\Gamma(n,3))=f(n,3)$.
\end{lemma}
\proof We shall use the notions in Construction~\ref{construction}. We claim that $v_d\in V(\Gamma(n,d))$. In fact, the claim holds if and only if $m\geq d$. If $k=1$, then $m=n-1\geq d$. If $k\geq 2$, then $(k-1)+d^{2(k-1)}<n$, which implies that
$m=n-k\geq d^{2k-2}\geq d^2>d.$
Hence, our claim is valid.

(i) Note that $k=f(n,d)=1$ and $V(\Gamma(n,d))=\{u_1\}\cup\{v_1,\ldots,v_m\}$ with
$$
v_j\in\{(a_1,a_2)\mid 1\leq a_1\leq d,1\leq a_2\leq d\},\qquad j=1,\ldots,m.
$$
By (F3), for distinct vertices $v_j=(a_1,a_2)$ and $v_l=(b_1,b_2)$, we have
$$
\partial(v_j,v_l)=\max\{|a_1-b_1|,|a_2-b_2|\}\leq d-1.
$$
Further more, by (F1) and (F2), one gets
$$
\begin{array}{rl}
&\partial(u_1,v_j)=\partial(u_1,(1,a_2)) +\partial((1,a_2),(a_1,a_2))=a_1,\\
&\partial(v_j,u_1)=\partial((a_1,a_2),(a_1,1)) +\partial((a_1,1),u_1)=a_2,
\end{array}
$$
and so $\tilde\partial(u_1,v_j)=(a_1,a_2)$. Hence, the set $\{u_1\}$ is a basis of $\Gamma(n,d)$. Noting that $v_d=(d,1)\in V(\Gamma)$, we obtain (i).

(ii) For any integers $a,b\in \{1,2\}$, we have $-1\leq a-b\leq 1$. Pick two distinct vertices $v_j=(a_1,a_2,\ldots,a_{2k})$ and $v_l=(b_1,b_2,\ldots,b_{2k})$.
It follows from (F3) that $\tilde\partial(v_j,v_l)=(1,1)$. Therefore, by (F1) and (F2), one gets $\tilde\partial(u_i,v_j)=(a_{2i-1},a_{2i})$, which implies that $\{u_1,\ldots,u_k\}$ is a weakly resolving set of $\Gamma(n,2)$, and so $\dim(\Gamma(n,2))=k=f(n,2)$ by Theorem~\ref{nd}.
Since $\tilde\partial(u_i,v_1)=(1,1)$ and $\partial(u_1,v_2)=2$, the diameter of $\Gamma(n,2)$ is $2$. Hence, (ii) follows.

(iii) Choose  distinct vertices $v_j=(a_1,a_2,\ldots,a_{2k})$ and $v_l=(b_1,b_2,\ldots,b_{2k})$. If there exists an index  $i\in\{1,2,\ldots,k\}$ such that $(a_{2i-1},b_{2i-1})=(1,3)$ or $(a_{2i},b_{2i})=(3,1)$, then $\partial(v_j,v_l)=2$; otherwise, we have $\partial(v_j,v_l)=1$. Therefore,  one gets $\tilde\partial(u_i,v_j)=(a_{2i-1},a_{2i})$, which implies that $\{u_1,\ldots,u_k\}$ is a weakly resolving set of $\Gamma(n,3)$, and so $\dim(\Gamma(n,3))=k=f(n,3)$ by Theorem~\ref{nd}. Note that $(u_i,v_1,u_j)$ is a path for any $i,j\in\{1,2,\ldots,k\}$. Since $\partial(u_1,v_3)=3$, the diameter of $\Gamma(n,2)$ is $3$, as desired.
$\qed$

Let $\mathcal{G}_{k,d}$ denote the set of $k$-dimensional digraphs with diameter $d$. For each $\Gamma\in\mathcal G_{k,d}$, by Theorem~\ref{nd}, we have sharp bounds
\begin{equation}\label{nkd}
k+d\leq|V(\Gamma)|\leq d^{2k}+k.
\end{equation}
In the rest of this section, we shall establish sharp bounds for $|A(\Gamma)|$.
For $x,y\in V(\Gamma)$, we say that $(x,y)$ and $(y,x)$ are {\em symmetric arcs} if $\tilde\partial(x,y)=(1,1)$.
Denote by $\overline\Gamma(d^{2k}+k,d)$ the digraph obtained from $\Gamma(d^{2k}+k,d)$ by adding symmetric arcs between $u_i$ and $u_j$ for $1\leq i<j\leq k$, where $u_i$ and $u_j$  refer to (\ref{uv}).

\begin{lemma}\label{embed}
Let $k$ and $d$ be positive integers. Pick a digraph $\Gamma\in\mathcal G_{k,d}$.

{\rm(i)} The digraph $\Gamma$ is isomorphic to a subdigraph of $\overline{\Gamma}(d^{2k}+k,d)$.

{\rm(ii)} If $k\geq 2$, $d=3$ and $|V(\Gamma)|=3^{2k}+k$, then $\Gamma$ is isomorphic to  a subdigraph of $\Gamma(3^{2k}+k,3)$.

{\rm(iii)} If $k\geq 2$ and $d\geq 4$, then $|V(\Gamma)|<d^{2k}+k$.
\end{lemma}
\proof
Let $X=\{x_1,\ldots,x_k\}$ be a basis of $\Gamma$ and $Y=V(\Gamma)\setminus\{x_1,\ldots,x_k\}$. For each vertex $y\in Y$, write
$$
\alpha_y=(\partial(x_1,y),\partial(y, x_1),\partial(x_2,y),\partial(y, x_2), \ldots, \partial(x_k,y),\partial(y, x_k)).
$$
Note that $1\leq \partial(x_i,y)\leq d$ and $1\leq\partial(y, x_i)\leq d$. Then $\alpha_y\in V(\overline\Gamma(d^{2k}+k,d))$.

(i)  Our goal is to establish a monomorphism from $\Gamma$ to $\overline{\Gamma}(d^{2k}+k,d)$. Define a map
$\iota: V(\Gamma)\longrightarrow V(\overline\Gamma(d^{2k}+k,d))$
as $\iota(x_i)=u_i$ and $\iota(y)=\alpha_y$, where $u_i$ refers to (\ref{uv}). This map is an injection because $X$ is a weakly resolving set.
 In order to show that $\iota$ is a homomorphism, we divide all arcs of $\Gamma$ in the following four cases.

{\em Case 1.} $(x_i,x_j)\in A(\Gamma)$. The fact is that $\overline{\Gamma}(d^{2k}+k,d)$ has an arc from $u_i$ to $u_j$.

{\em Case 2.} $(x_i,y)\in A(\Gamma)$. By (F1), the fact that $\partial(x_i,y)=1$ implies that $(u_i,\alpha_y)$ is an arc of $\Gamma(d^{2k}+k,d)$, as well as $\overline\Gamma(d^{2k}+k,d)$.

{\em Case 3.} $(y, x_i)\in A(\Gamma)$. By (F2), the fact that $\partial(y, x_i)=1$ implies that $(\alpha_y, u_i)$ is an arc of $\Gamma(d^{2k}+k,d)$, as well as $\overline\Gamma(d^{2k}+k,d)$.

{\em Case 4.} $(y, z)\in A(\Gamma)$ for $z\in Y$.
For each $x_i$, we have $\partial(x_i,y)-\partial(x_i,z)\geq -1$ and $\partial(y,x_i)-\partial(z,x_i)\leq 1$ by the triangle inequality. It follows from (F3) that $(\alpha_y,\alpha_z)$ is an arc of $\Gamma(d^{2k}+k,d)$, as well as $\overline\Gamma(d^{2k}+k,d)$.

Consequently, the digraph $\Gamma$ is isomorphic to a subdigraph of $\overline{\Gamma}(d^{2k}+k,d)$ since $\iota$ is a monomorphism from $\Gamma$ to $\overline{\Gamma}(d^{2k}+k,d)$.

(ii) We claim that $(x_i,x_j)\not\in A(\Gamma)$ for any  $x_i,x_j\in X$.
In fact, since $|Y|=3^{2k}$, we infer that
$\alpha_y$ ranges over the set
$$
\{(a_1,a_2,\ldots,a_{2k-1},a_{2k})\mid 1\leq a_r\leq 3\text{ for }r=1,2,\ldots,2k\}
$$
as $y$ ranging over all vertices in $Y$. Hence, for any two distinct vertices $x_i$ and $x_j$ in $X$, there exists a vertex $y\in Y$ such that $\partial (x_i,y)=3$ and $\partial(x_j,y)=1$, which implies that $(x_i,x_j)\not\in A(\Gamma)$ by the triangle inequality. Therefore, the claim is valid. An argument similar to the one used in the proof of (i), except {\em Case 1},  shows that there exists a monomorphism from $\Gamma$ to $\Gamma(3^{2k}+k,3)$, and so (ii) holds.

(iii) By contradiction,
suppose that $|V(\Gamma)|\geq d^{2k}+k$.
By (\ref{nkd}), we have $|Y|=|V(\Gamma)|-k=d^{2k}$.
Thus, as $y$ ranging over all vertices in $Y$, the vector $\alpha_y$ ranges over the set
$$
\{(a_1,a_2,\ldots,a_{2k-1},a_{2k})\mid 1\leq a_r\leq d\text{ for } r=1,2,\ldots,2k\}.
$$
Therefore, there exist vertices $y_1,y_2\in Y$ such that
$$
\partial(x_1,y_1)=4,\;\partial(x_2,y_1)=1 \text{ and } \partial(x_1,y_2)=\partial(y_2,x_2)=1,
$$
which implies that
$$
4=\partial(x_1,y_1)\leq \partial(x_1,y_2)+\partial(y_2,x_2)+\partial(x_2,y_1)=3,
$$
a contradiction.
$\qed$

As refer to Construction~\ref{construction}, let $e(k,d)$ denote the number of arcs in family (F3) in the digraph $\Gamma(d^{2k}+k,d)$.

\begin{lemma}\label{ekd}
$e(k,d)=(\frac{d^2+3d-2}{2})^{2k}-d^{2k}$.
\end{lemma}
\proof As refer to (\ref{uv}), let $\Delta$ denote the induced subdigraph of $\Gamma(d^{2k}+k,d)$ on the  vertex subset $\{v_1,\ldots,v_{d^{2k}}\}$. Then
$$
e(k,d)=|A(\Gamma(d^{2k}+k,d))|=\sum_{v_j\in V(\Delta)}\delta^+(v_j),
$$
where $\delta^+(v_j)$ is the out-degree of $v_j$  in the digraph $\Delta$. Note that
$$
\delta^+((a_1,a_2,\ldots,a_{2k-1},a_{2k}))=\prod_{r=1}^{2k} \theta(a_r)-1,
$$
where
$$
\theta(a_r)=\left\{
\begin{array}{ll}
\min\{a_{r}+1,d\}, & \text{if $r$ is odd},\\
\min\{d-a_{r}+2,d\}),& \text{if $r$ is even}.
\end{array}\right.
$$
Write $\mu=\sum_{a_r=1}^d\theta(a_r)=2+3+\cdots+d+d=\frac{d^2+3d-2}{2}$.
Hence,
\begin{eqnarray*}
\sum_{v_j\in V(\Delta)}\delta^+(v_j)&=& \sum_{a_1,\ldots,a_{2k-1}=1}^d(\sum_{a_{2k}=1}^d \delta^+((a_1,\ldots,a_{2k-1},a_{2k}))) \\
&=&\sum_{a_1,\ldots,a_{2k-2}=1}^d(\sum_{a_{2k-1}=1}^d (\prod_{r=1}^{2k-1}\theta(a_r)\mu-d)) \\ &=&\sum_{a_1,\ldots,a_{2k-3}=1}^d(\sum_{a_{2k-2}=1}^d (\prod_{r=1}^{2k-2}\theta(a_r)\mu^2-d^2))\\
&=&\cdots\cdots=\mu^{2k}-d^{2k},
\end{eqnarray*}
as desired.
$\qed$

\begin{thm}\label{ag}
  Let $\Gamma$ be a $k$-dimensional digraph with diameter $d$. Then
  \begin{equation}\label{akd}
    k+d\leq |A(\Gamma)|\leq (\frac{d^2+3d-2}{2})^{2k}+(2k-d)d^{2k-1}+k^2-k.
  \end{equation}
  Moreover,

  {\rm(i)} the lower bound is attained if and only if $k=1$ and $\Gamma$ is isomorphic to a directed cycle with length $d+1$.

  {\rm(ii)} the upper bound is attained if and only if one of the following holds.

  \hspace{5mm}{\rm(a)} $k=1$ and $\Gamma$ is isomorphic to $\Gamma(d^2+1,d)$.

  \hspace{5mm}{\rm(b)} $d=1$ and $\Gamma$ is isomorphic to the complete digraph of order $k+1$.

  \hspace{5mm}{\rm(c)} $d=2$ and $\Gamma$ is isomorphic to $\overline\Gamma(4^k+k,2)$.
\end{thm}
\proof A short calculation reveals
$$
|A(\overline\Gamma(d^{2k}+k,d))|=k(k-1)+2kd^{2k-1}+e(k,d),
$$
which is equal to the upper bound in (\ref{akd}) by Lemma~\ref{ekd}.
Note that $|V(\Gamma)|\leq|A(\Gamma)|$. According to the inequality  (\ref{nkd}) and Lemma~\ref{embed}, we get the inequality (\ref{akd}).

(i) Observe that $|V(\Gamma)|=|A(\Gamma)|$ if and only if $\Gamma$ is a directed cycle.
We obtain (i) since a directed cycle is $1$-dimensional.

(ii) Note that $\Gamma(d^2+1,d)$ is the same as $\overline\Gamma(d^2+1,d)$. It follows from Lemma~\ref{eg2}~(i) that $\overline\Gamma(d^2+1,d)\in\mathcal G_{1,d}$.
It is clear that the complete digraph of order $k+1$ is isomorphic to $\overline\Gamma(1+ k,1)$ and $\overline\Gamma(1+ k,1)\in\mathcal G_{k,1}$.
An argument similar to the one used in the proof of Lemma~\ref{eg2}~(ii) shows that
$\overline\Gamma(4^k+k,2)\in\mathcal G_{k,2}$.
Hence, the ``if'' implication follows, while the ``only if'' implication follows immediately from Lemma~\ref{embed}.
$\qed$

\section{$1$-dimensional digraphs}
Recall that $\mathcal G_{1,d}$ denotes the set of $1$-dimensional digraphs with diameter $d$. It follows from Lemma~\ref{eg2}~(i) that $\mathcal G_{1,d}$ is not empty for each positive integer $d$. We list some properties of digraphs $\Gamma\in\mathcal G_{1,d}$ that have been stated in Lemma~\ref{embed}, Theorems~\ref{nd} and~\ref{ag}.

\begin{prop}\label{dim1}
  Let $\Gamma\in\mathcal G_{1,d}$.

  {\rm(i)} The digraph $\Gamma$ is isomorphic to a subdigraph of $\Gamma(d^2+1,d)$.

  {\rm(ii)} $d+1\leq|V(\Gamma)|\leq d^2+1$.

  {\rm(iii)} $d+1\leq |A(\Gamma)|\leq (\frac{d^2+3d-2}{2})^2+d(2-d)$.
\end{prop}

All bounds in Proposition~\ref{dim1} are sharp, and in particular, the equalities in (iii) have been characterized in Theorem~\ref{ag} for $k=1$.
Example~\ref{eg1} for $n-d=1$, and Lemma~\ref{eg2}~(i) give pairwise non-isomorphic digraphs with weak metric dimension $1$. The following example shows that, non-isomorphic digraphs with the same vertices and the same diameter may have a common basis even if they have the same girth, that is the length of a shortest directed cycle in this digraph.

\begin{figure}
  \centering
  % Requires \usepackage{graphicx}
  \includegraphics[width=12cm]{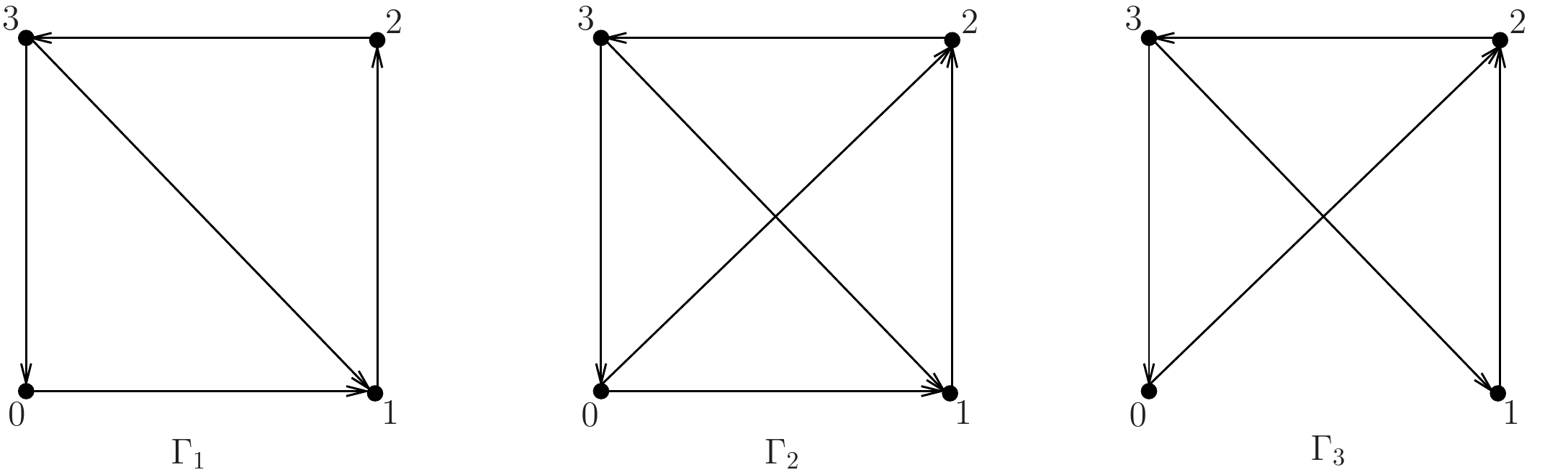}\\
  \caption{Non-isomorphic $1$-dimensional graphs of order $4$, diameter $3$ and girth $3$}\label{fig1}
\end{figure}

\begin{eg}
Let $\Gamma$ be the directed cycle with the vertex set $\{0,1,2,3\}$ and the arc set $\{(0,1),(1,2),(2,3),(3,0)\}$.
As refer to Figure~\ref{fig1}, $\Gamma_1$ is obtained from $\Gamma$ by adding the arc $(3,1)$; $\Gamma_2$ is obtained from $\Gamma_1$ by adding the arc $(0,2)$; $\Gamma_3$ is obtained from $\Gamma_2$ by deleting the arc $(0,1)$. All these  digraphs have the same vertices, the same diameter and a common basis $\{0\}$. All of digraphs  $\Gamma_1$, $\Gamma_2$ and $\Gamma_3$ have girth $3$. Digraphs $\Gamma_1$ and $\Gamma_3$ are not isomorphic because $\Gamma_1$ has a directed cycle with length $4$ and $\Gamma_3$ does not have.
\end{eg}

It seems very difficult to characterize all $1$-dimensional digraphs. In the rest of  this section, we characterize all $1$-dimensional  vertex-transitive  digraphs.

Denote by Aut$(\Gamma)$ the full automorphism group of a digraph $\Gamma$. We say that $\Gamma$ is {\em vertex-transitive} if Aut$(\Gamma)$ acts transitively on the vertex set $V(\Gamma)$.
Cayley digraphs are a family of classical vertex-transitive digraphs.
Let $G$ be a finite group and let $S$ be a generating subset not containing
the identity element. The {\em Cayley digraph} Cay$(G,S)$ is a digraph with the vertex set $G$ and the arc set $\{(x,xs)\mid x\in G, s\in S\}$.
Let $\mathbb{Z}_n$ denote the additive group of integers modulo $n$.
We now state the main theorem of this section, which will be proved at the end of this section.

\begin{thm}\label{vertextran1dim}
  A vertex-transitive digraph $\Gamma$ of order at least $2$ is $1$-dimensional if and only if $\Gamma$ is isomorphic to one of the following digraphs:

   {\rm(i)} {\rm Cay}$(\mathbb{Z}_{n}, \{1\})$, where $n\geq 2$.

 {\rm(ii)} {\rm Cay}$(\mathbb{Z}_{2n}, \{1,2\})$, where $n\geq 2$.

 {\rm(iii)} {\rm Cay}$(\mathbb{Z}_{2}\oplus \mathbb{Z}_{n},\{(1,0),(0,1)\})$, where $n\geq 3$.

 {\rm(iv)} {\rm Cay}$(\mathbb{Z}_{2}\oplus \mathbb{Z}_{2n},\{(1,0),(0,1),(0,2)\})$, where $n\geq 3$.
\end{thm}

Actually, digraphs in the above theorem are as a classification of so-called thin weakly distance-regular digraphs (see Proposition~\ref{wdr}). Our proof of  Theorem~\ref{vertextran1dim} depends on  this result.
To prove Theorem~\ref{vertextran1dim}, we begin by giving notations. For a digraph $\Gamma$, write $\tilde\partial(\Gamma)=\{ \tilde\partial (x,y)\mid x,y\in V(\Gamma)\}$ and let
$$
\Gamma_{\tilde i}(x)=\{y\in V(\Gamma)\mid \tilde\partial(x,y)=\tilde i\}
$$
for $\tilde i\in\tilde\partial(\Gamma)$ and $x\in V(\Gamma)$.

\begin{ob}\label{ob1}
  Let $x$ be a vertex of a digraph $\Gamma$. Then $\{x\}$ is a weak metric basis of $\Gamma$ if and only if  $|\Gamma_{\tilde i}(x)|\leq 1$ for each $\tilde i\in\tilde\partial(\Gamma)$.
\end{ob}

A digraph $\Gamma$ is \emph{weakly distance-regular} if, for any $\tilde{h}$, $\tilde{i}$, $\tilde{j}\in\tilde{\partial}(\Gamma)$, the number of $z\in V(\Gamma)$ such that $\tilde{\partial}(x,z)=\tilde{i}$ and $\tilde{\partial}(z,y)=\tilde{j}$ is constant whenever $\tilde{\partial}(x,y)=\tilde{h}$. The number $|\Gamma_{\tilde i}(x)|$ does not depend on the choice $x\in V(\Gamma)$, which is denoted by $k_{\tilde i}$. A weakly distance-regular digraph $\Gamma$ is {\em thin} if $k_{\tilde i}=1$ for any $\tilde{i}\in\tilde\partial(\Gamma)$. Wang and Suzuki~\cite{wangsuzuki} first introduced the concept of weakly distance-regular digraphs. For more results on weakly distance-regular digraphs, see \cite{suzuki,KSW04,YYF16,YYF18,YYF20,YYF20+}.

Now we state the result of classifying all thin weakly distance-regular digraphs.

\begin{prop}{\rm\cite[Theorem 1.2]{suzuki}}\label{wdr}
A digraph $\Gamma$ is a thin weakly distance-regular digraph of order at least $2$ if and only if $\Gamma$ is isomorphic to one of  digraphs in Theorem~\ref{vertextran1dim}.
\end{prop}

A digraph $\Gamma$ is  {\em weakly distance-transitive}  if, for any vertices
$x,y,x'$ and $y'$ with $\tilde\partial(x,y)=\tilde\partial(x',y')$, there exists
an automorphism $\sigma\in {\rm Aut}(\Gamma)$ such that $\sigma(x)=x'$ and $\sigma(y)=y'$.

\begin{ob}\label{ob2}
  A   weakly distance-transitive  digraph is weakly distance-regular.
\end{ob}

\begin{lemma}\label{lemma1}
  If $\Gamma$ is a $1$-dimensional vertex-transitive digraph, then $\Gamma$ is weakly distance-transitive.
\end{lemma}
\proof Pick any vertices $x,y,x'$ and $y'$ with $\tilde\partial(x,y)=\tilde\partial(x',y')$. Note that there exists an automorphism $\sigma$ of $\Gamma$ such that $\sigma(x)=x'$. Then
$$
\tilde\partial(x',y')=\tilde\partial(x,y)=\tilde\partial(x',\sigma(y)).
$$
Since $\{x'\}$ is a weakly resolving set of $\Gamma$, one has $y'=\sigma(y)$, as desired.
$\qed$

{\noindent\em Proof of Theorem~\ref{vertextran1dim}:}  According to Lemma~\ref{lemma1}, Observations~\ref{ob2},~\ref{ob1} and Proposition~\ref{wdr}, the ``only if'' implication follows, while the ``if'' implication follows from Proposition~\ref{wdr}, Observation~\ref{ob1} and the fact that Cayley digraphs are vertex-transitive.
$\qed$

\section{$(n-2)$-dimensional digraphs of order $n$}

We have characterized  $(n-1)$-dimensional digraphs of order $n$ in Corollary~\ref{complete}. In fact, they are precisely the complete digraphs.  In this section, we classify all $(n-2)$-dimensional digraphs of order $n$.

An arc $(x,y)$ is of {\em type} $(1,r)$ if $\partial(y,x)=r$. A digraph is said to be an {\em undirected graph}, or a {\em graph} for simplify, if each arc is of type $(1,1)$. In other words, each arc in a graph has a symmetric arc. Let $K_t$ and $\overline K_t$ denote the complete and null digraph of order $t$, respectively.
 For  digraphs $G$ and $H$, we use $G\cup H$ to denote the disjoined union of $G$ and $H$, and $G+H$ to denote the digraph obtained from $G\cup H$ by adding symmetric arcs between each vertex of $G$ and each vertex of $H$. Note that the weak metric dimension of a graph is equal to the metric dimension of this graph.
Graphs of order $n$ having metric dimension $n-2$ have been classified.

\begin{prop}{\rm\cite[Theorem 4]{chartrandejo}}\label{graph}
  Let $\Gamma$ be a connected graph of order $n\geq 4$. Then $\dim(\Gamma)=n-2$ if and
only if $\Gamma$ is isomorphic to $\overline K_t+\overline K_{n-t}$, $K_t+\overline{K}_{n-t}$ or $K_t+(K_1\cup K_{n-t-1})$, where $1\leq t\leq n-2$.
\end{prop}

Let $G$ be a digraph with the vertex set $\{0,1,\ldots,m-1\}$, and let $H_0,H_1,\ldots,H_{m-1}$ be digraphs.  The {\em generalized lexicographic product} of $H_0,H_1,\ldots,H_{m-1}$ by $G$, denoted by $G[H_0,H_1,\ldots,H_{m-1}]$, is the digraph with the vertex set $V(H_0)\cup V(H_1)\cup\cdots\cup V(H_{m-1})$, where there is an arc from a vertex $x\in H_i$ to a vertex $y\in H_j$ if and only if one of the following conditions holds:

(C1) If $i=j$, then $(x,y)$ is an arc of $H_i$.

(C2) If $i\neq j$, then $(i,j)$ is an arc of $G$.

\medskip

Now we state the main theorem of this section.

\begin{figure}
  \centering
  % Requires \usepackage{graphicx}
  \includegraphics[width=13.5cm]{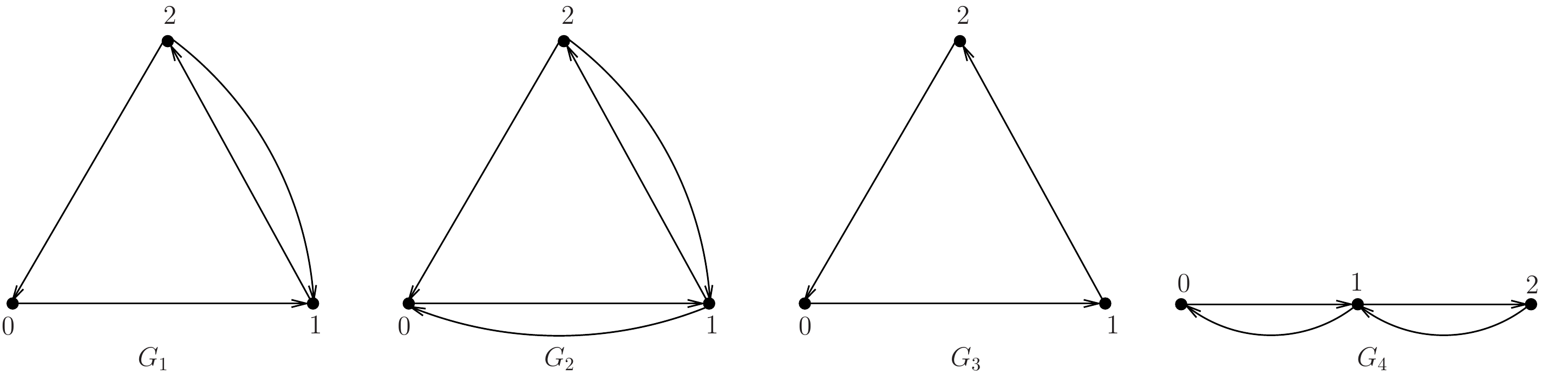}\\
  \caption{All $1$-dimensional digraphs of order $3$}\label{3tu}
\end{figure}

\begin{thm}\label{n-2main}
  Let $\Gamma$ be a digraph of order $n\geq 4$. Suppose that $\Gamma$ is not an undirected graph. Then $\dim(\Gamma)=n-2$ if and only if $\Gamma$ is isomorphic to $G_2[K_1,P_2,K_1]$, $G_1[K_1,K_t,K_{n-t-1}]$, $G_2[K_1,K_t,K_{n-t-1}]$, $G_2[K_t,K_{n-t-1},K_1]$ or $G_2[K_1,\overline K_{n-2},K_1]$ with $1\leq t\leq n-2$, where $P_2$ is a directed path with length $1$ and $G_1$, $G_2$ are as refer to Figure~\ref{3tu}.
\end{thm}

\noindent {\bf Remark} Any $1$-dimensional digraph of order $3$ is isomorphic to $G_1$, $G_2$, $G_3$ or $G_4$ in Figure~\ref{3tu}.

\medskip

In the rest of this section, we shall prove Theorem~\ref{n-2main}.

\begin{lemma}\label{gik1ktks}
For $i=1,2$, we have $\dim(G_i[K_1,K_t,K_s])=s+t-1$.
\end{lemma}
\proof Let $V(K_1)=\{u\}$, $V(K_t)=\{v_1,\ldots,v_t\}$ and $V(K_s)=\{w_1,\ldots,w_s\}$.
Distinct vertices $v_i$ and $v_j$ can only be resolved by $v_i$ or $v_j$, and distinct vertices $w_i$ and $w_j$ also can only be resolved by $w_i$ or $w_j$. Note that $G_i$ has symmetric arcs between $1$ and $2$. Then $v_i$ and $w_j$ can only be resolved by $u$, $v_i$ or $w_j$. To sum up, we have derived that $\dim(G_i[K_1,K_t,K_s])\geq t-1+s-1+1=s+t-1$. Since $G_i[K_1,K_t,K_s]$ is not complete and has $s+t+1$ vertices, we get the desired result.
$\qed$

Note that $G_2[K_t,K_{n-t-1},K_1]$ is obtained from $G_2[K_1,K_{n-t-1},K_t]$  by reversing directions of all arcs. Hence, $\dim(G_2[K_t,K_{n-t-1},K_1])=\dim(G_2[K_1,K_{n-t-1},K_t])$. It is routine to check $\dim(G_2[K_1, P_2,K_1])=2$ and
$\dim(G_2[K_1,\overline K_{n-2},K_1])=n-2$. The following result follows from Lemma~\ref{gik1ktks}.

\begin{lemma}\label{check}
  All digraphs in Theorem~\ref{n-2main} are $(n-2)$-dimensional.
\end{lemma}

\begin{lemma}\label{type}
  If a digraph $\Gamma$ has an arc of type $(1,r)$, then $\dim(\Gamma)\leq |V(\Gamma)|-r$.
\end{lemma}
\proof Let $(x,y)$ be an arc of type $(1,r)$. Then there exists a shortest path $(y=y_0,y_1,\ldots,y_r=x)$. Note that $\partial (y_0,y_i)=i$ for any $i\in\{1,\ldots,r\}$. Hence, for distinct vertices $y_i$ and $y_j$ in $\{y_1,\ldots,y_r\}$, the vertex $y_0$ resolves $y_i$ and $y_j$, which implies that   $V(\Gamma)\setminus\{y_1,\ldots,y_r\}$ is a weakly resolving set of $\Gamma$, and so $\dim(\Gamma)\leq |V(\Gamma)|-r$.
$\qed$

The following result is immediate from Lemma~\ref{type}.

\begin{lemma}\label{12}
  If $\dim(\Gamma)=|V(\Gamma)|-2$, then each arc is of type $(1,1)$ or $(1,2)$.
\end{lemma}

To prove the necessity in Theorem~\ref{n-2main}, on account of Lemma~\ref{12}, we always suppose that $\Gamma$ is an $(n-2)$-dimensional digraph of order $n\geq 4$, and has an arc $(x_2,x_0)$ of type $(1,2)$ and a directed path $(x_0,x_1,x_2)$. Write
$$
X_1=\{x\in V(\Gamma)\mid \partial(x_0,x)=1\}\quad\text{ and }\quad X_2=\{x\in V(\Gamma)\mid \partial (x_0,x)=2\}.
$$
Then both $X_1$ and $X_2$ are nonempty since $x_1\in X_1$ and $x_2\in X_2$. By the upper bound in Theorem~\ref{nd}, the following result is valid.

\begin{lemma}\label{diam}
 The diameter of $\Gamma$ is $2$ and $V(\Gamma)=\{x_0\}\cup X_1\cup X_2$.
\end{lemma}

\begin{lemma}\label{3vertices}
 Pick three pairwise distinct vertices $x,y,z\in V(\Gamma)$.

  {\rm(i)} For any vertex $w\in V(\Gamma)\setminus\{x,y,z\}$, the two way distances $\tilde\partial(w,x),\; \tilde\partial(w,y)$ and $\tilde\partial(w,z)$
  are not pairwise distinct.

  {\rm (ii)}  If $\{x,y,z\}\subseteq X_1$, then $\tilde\partial(x,y)=\tilde\partial(x,z)=(1,1)$ or $(2,2)$.

  {\rm (iii)}If $\{x,y,z\}\subseteq X_2$, then $\tilde\partial(x,y)=\tilde\partial(x,z)=(1,1)$.
\end{lemma}
\proof (i) By contradiction, suppose that $\tilde\partial(w,x),\; \tilde\partial(w,y)$ and $\tilde\partial(w,z)$ are pairwise distinct. Then $V(\Gamma)\setminus\{x,y,z\}$ is a weakly resolving set of $\Gamma$, which contradicts the fact that  $\dim(\Gamma)=n-2$.

(ii) Note that $x_2\in X_2$. If $\tilde\partial(x,y)\neq \tilde\partial(x,z)$, then $V(\Gamma)\setminus\{y,z,x_2\}$ is a weakly resolving set of $\Gamma$ since $\partial(x_0,y)=\partial(x_0,z)\neq \partial(x_0,x_2)$. This contradicts the fact that  $\dim(\Gamma)=n-2$. Thus, one has $\tilde\partial(x,y)=\tilde\partial(x,z)$.
Using the same argument, we get $\tilde\partial(y,x)=\tilde\partial(y,z)$ and $\tilde\partial(z,x)=\tilde\partial(z,y)$. Therefore, we infer that $\tilde\partial(x,y)=\tilde\partial(y,x)$, which implies that $\partial(x,y)=\partial(y,x)$, and so $\tilde\partial(x,y)=(1,1)$ or $(2,2)$ by Lemma~\ref{diam}.

(iii) The proof of
$$
\tilde\partial(x,y)=\tilde\partial(x,z)=(1,1)\text{ or }(2,2)
$$
 is almost identical with the proof of (ii), the major change being the substitution of $x_2$ by $x_1$. Suppose $\tilde\partial(x,y)=(2,2)$. Then the two way distance from any vertex in $X_2$ to another vertex in $X_2$ is $(2,2)$. Note that $|X_2|\geq 3$. Pick $w\in X_2\setminus\{x_2\}$. Then
$\tilde\partial(x_2,x_0),\;\tilde\partial(x_2,x_1)$ and $\tilde\partial(x_2,w)$
are pairwise distinct, contrary to (i). Hence, we get $\tilde\partial(x,y)=\tilde\partial(x,z)=(1,1)$.
$\qed$

\begin{lemma}\label{x1x2}
  {\rm(i)} For each $w\in X_1$, arcs $(x_0,x_1)$ and $(x_0,w)$ have the same type.

  {\rm(ii)} For each $w\in X_2$, there is an arc of type $(1,2)$ from $w$ to $x_0$ in $\Gamma$.
\end{lemma}
\proof (i) Only need to consider $w\neq x_1$. In this case, vertices $x_0,x_1,w,x_2$ are pairwise distinct. It follows from Lemma~\ref{3vertices}~(i) that $\tilde\partial(x_0,x_1)$, $\tilde\partial(x_0,w)$ and $\tilde\partial(x_0,x_2)$ are not pairwise distinct, which implies that $\tilde\partial(x_0,x_1)=\tilde\partial(x_0,w)$, and so (i) holds.

(ii) An argument similar to the one used in the proof of (i) shows that $\tilde\partial(x_0,x_2)=\tilde\partial(x_0,w)$, and so $\tilde\partial(w,x_0)=\tilde\partial(x_2,x_0)=(1,2)$, as desired.
$\qed$

A {\em clique} in a digraph is a vertex subset $S$ such that the two way distance from each vertex in $S$ to another vertex in $S$ is $(1,1)$.
An {\em independent set} in a digraph is a vertex subset $R$ such that there are no arcs from any vertex in $R$ to any vertex in $R$. The following result is immediate from Lemma~\ref{3vertices}~(ii) and (iii).

\begin{lemma}\label{cliqueindependent}
{\rm(i)} If $|X_1|\geq 3$, then $X_1$ is a clique or an independent set.

{\rm(ii)} If $|X_2|\geq 3$, then $X_2$ is a clique.
\end{lemma}

\begin{prop}\label{case1}
  If $X_1$ is an independent set with at least $2$ vertices, then $\Gamma$ is isomorphic to $G_2[K_1,\overline K_{n-2},K_1]$.
\end{prop}
\proof Take any $w\in X_1\setminus\{x_1\}$. Then $\tilde\partial(x_1,w)=(2,2)$ by Lemma~\ref{diam}. In the following, we divide the proof in four steps.

Our first goal is to show that
\begin{equation}\label{4}
\tilde\partial(x_0,x_1)=\tilde\partial(x_0,w)= \tilde\partial(x_1,x_2)=(1,1).
\end{equation}
It follows from Lemma~\ref{3vertices}~(i)  that $\tilde\partial(x_1,x_0)$, $\tilde\partial(x_1,w)$ and $\tilde\partial(x_1,x_2)$ are not pairwise distinct,
which implies that $\tilde\partial(x_1,x_0)=\tilde\partial(x_1,x_2)=(1,1)$, and so $\tilde\partial(w,x_0)=(1,1)$ by Lemma~\ref{x1x2}~(i). Hence, we obtain (\ref{4}).

The next thing to do is to prove $X_2=\{x_2\}$.
Considering $\tilde\partial(w,x_0)$, $\tilde\partial(w,x_1)$ and $\tilde\partial(w,x_2)$,
by Lemma~\ref{3vertices}~(i), we have
\begin{equation}\label{5}
  \tilde\partial(w,x_2)=(1,1) \text{ or }(2,2).
\end{equation}
Suppose that $X_2\neq\{x_2\}$. Take $y\in X_2\setminus\{x_2\}$.
Then $V(\Gamma)\setminus\{x_0,w,x_2\}$ contains $\{x_1,y\}$ and is not a weakly resolving set of $\Gamma$ since $\dim(\Gamma)=n-2$. It follows from (\ref{4}) that $x_1$ resolves both $\{x_0,w\}$ and $\{w,x_2\}$,
which implies that $y$ can not resolve $x_2$ and $x_0$, and so
$\tilde\partial(y,x_2)=\tilde\partial(y,x_0)=(1,2)$ from Lemma~\ref{x1x2}~(ii).
Therefore, we have derived from (\ref{5}) that $\tilde\partial(x_2,w)$, $\tilde\partial(x_2,y)$ and $\tilde\partial(x_2,x_0)$ are pairwise distinct, contrary to Lemma~\ref{3vertices}~(i).
Hence, we get $X_2=\{x_2\}$.

Another step is to prove $\tilde\partial(w,x_2)=(1,1)$. Suppose for the contrary that
$\tilde\partial(w,x_2)\neq (1,1)$. By (\ref{5}), we get $\tilde\partial(w,x_2)=(2,2)$.
Then there is a vertex $z\in V(\Gamma)$ such that $(w,z,x_2)$ is a directed path.
Since $X_2=\{x_2\}$ and $\partial(x_0,x_2)=2$, one has $z\in X_1$, contrary to the condition that $X_1$ is independent.  Thus,  $\tilde\partial(w,x_2)=(1,1)$.

Finally, one obtains $|X_1|=n-2$ owing to $|X_2|=1$.
In view of (\ref{4}) and $\tilde\partial(w,x_2)=(1,1)$,   we conclude that $\Gamma$ is isomorphic to $G_2[K_1,\overline K_{n-2},K_1]$.
$\qed$

\begin{prop}\label{case2}
  Suppose that $X_1$ is a clique and $|X_2|=1$.

  {\rm (i)} If the arc $(x_0,x_1)$ is of type $(1,1)$, then $\Gamma$ is isomorphic to
  $$
  G_1[K_{1},K_t,K_{n-t-1}]\qquad \text{or}\qquad G_2[K_{t},K_{n-t-1},K_1],
  $$
   where $1\leq t\leq n-2$.

  {\rm (ii)} If the arc $(x_0,x_1)$ is of type $(1,2)$, then $\Gamma$ is isomorphic to $G_1[K_{1},K_{n-2},K_1]$.
\end{prop}
\proof Noting that $n\geq 4$, we have $|X_1|=n-2\geq 2$. For any $z\in X_1\setminus\{x_1\}$, we infer from Lemma~\ref{3vertices}~(i) that $\tilde\partial(x_2,x_0)$, $\tilde\partial(x_2,x_1)$ and $\tilde\partial(x_2,z)$ are not pairwise distinct. Then $\tilde\partial(x_2,z)=\tilde\partial(x_2,x_0)$ or $\tilde\partial(x_2,x_1)$. Write
$$
Y_1=\{y\in V(\Gamma)\mid \tilde\partial(x_2,y)=\tilde\partial(x_2,x_0)\}
\text{ and }
Y_2=\{y\in V(\Gamma)\mid \tilde\partial(x_2,y)=\tilde\partial(x_2,x_1)\}.
$$
Let $t=|Y_1|$. Since $Y_1$ and $Y_2$ are nonempty, we have $1\leq t\leq n-2$.

(i) Note that $Y_1\cup Y_2=X_1\cup\{x_0\}$, which is a clique on account of Lemma~\ref{x1x2}~(i).
If the arc $(x_1,x_2)$ is of type $(1,1)$, then $\Gamma$ is isomorphic to $G_2[K_{t},K_{n-t-1},K_1]$.
Suppose that the arc $(x_1,x_2)$ is of type $(1,2)$. It follows that $\Gamma$ is isomorphic to
$H[K_{t},K_{n-t-1},K_1]$, where $H$ is a digraph with the vertex set $\{0,1,2\}$ and the arc set $\{(0,1),(1,0),(1,2),(2,0)\}$. Since $H$ is isomorphic to $G_1$, we conclude that $\Gamma$ is isomorphic to $G_1[K_{1},K_t,K_{n-t-1}]$.

(ii) Since $\tilde\partial(x_0,x_1)=(1,2)$, the two way distance from $x_0$ to each vertex of $X_1$ is $(1,2)$ by Lemma~\ref{x1x2}~(i). Then $x_2$ is the unique vertex at distance $1$ to $x_0$ in $\Gamma$, which implies that for each $z\in X_1$, there exists an arc from $z$ to $x_2$. Choose $z_1\in X_1\setminus\{z\}$. By Lemma~\ref{3vertices}~(i), the two way distance $\tilde\partial(z,z_1)$, $\tilde\partial(z,x_0)$ and $\tilde\partial(z,x_2)$ are not pairwise distinct, and so $\tilde\partial(z,x_2)=(1,1)$. Hence $X_1\cup\{x_2\}$ is a clique, and then $\Gamma$ is isomorphic to $G_1[K_{1},K_{n-2},K_1]$.
$\qed$

\begin{prop}\label{case3}
Suppose that $X_1$ is neither a clique nor an independent set. If $|X_2|=1$, then $\Gamma$ is isomorphic to $G_2[K_1,P_2,K_1]$.
\end{prop}
\proof It follows from Lemma~\ref{cliqueindependent}~(i) that $|X_1|=2$. Then $\Gamma$ has exactly $4$ vertices. It is straightforward to get the desired result.
$\qed$

\begin{lemma}\label{arcx1x2}
 If $|X_1|\geq 2$ and $|X_2|\geq 2$, then $\tilde\partial(w_1,w_2)=\tilde\partial(x_1,x_2)$ for each $w_1\in X_1$ and each $w_2\in X_2$.
\end{lemma}
\proof We claim that $\tilde\partial(x_1,w_2)=\tilde\partial(x_1,x_2)$. In fact,  we only need to consider $w_2\neq x_2$. Pick a vertex $y_1\in X_1\setminus\{x_1\}$.
Note that $V(\Gamma)\setminus\{y_1,x_2,w_2\}$ contains $\{x_0,x_1\}$ and is not a weakly resolving set of $\Gamma$.
Since $\partial(x_0,y_1)=1\neq 2=\partial(x_0,x_2)=\partial(x_0,w_2)$, the vertex $x_1$ can not resolve $x_2$ and $w_2$. Hence, the claim is valid.

Now suppose that $\tilde\partial(w_1,w_2)\neq \tilde\partial(x_1,x_2)$. By the claim, we have $\tilde\partial(w_2,x_1)=\tilde\partial(x_2,x_1)\neq\tilde\partial(w_2,w_1)$.
Take a vertex $y_2\in X_2\setminus\{w_2\}$. Since $\partial(x_0,x_1)=\partial(x_0,w_1)=1\neq 2=\partial(x_0,y_2)$ and $w_2$ resolves $\{x_1,w_1\}$, we conclude that $V(\Gamma)\setminus\{x_1,w_1,y_2\}$ is a weakly resolving set of $\Gamma$, a contradiction.
Consequently, one gets $\tilde\partial(w_1,w_2)=\tilde\partial(x_1,x_2)$.
$\qed$

\begin{lemma}\label{x2clique}
 The set $X_2$ is a clique.
\end{lemma}
\proof By Lemma~\ref{cliqueindependent}, we only need to consider $|X_2|=2$. If $|X_1|=1$, it is routine to verify that $X_2$ is a clique. In the following, suppose $|X_1|\geq 2$.
Observe that for $i=1,2$, there exist vertices $y_i,z_i\in X_i$ such that $\tilde\partial(y_i,z_i)\neq(1,2)$.

By Lemma~\ref{arcx1x2}, one gets $\tilde\partial(y_2,z_1)=(\partial(x_2,x_1),1)$.
Noting that $\tilde\partial(y_2,x_0)=(1,2)\neq\tilde\partial(y_2,z_2)$, we have $\tilde\partial(y_2,z_2)=\tilde\partial(y_2,z_1)$ by Lemma~\ref{3vertices}~(i).
Suppose that $\tilde\partial(y_2,z_2)\neq (1,1)$. It follows that $\tilde\partial(y_2,z_2)=\tilde\partial(y_2,z_1)=(2,1)$, which implies that $y_2$ resolves both $\{x_0,z_1\}$ and $\{x_0,z_2\}$. Because of  the fact that $\{y_1,y_2\}\subseteq V(\Gamma)\setminus\{x_0,z_1,z_2\}$ and $V(\Gamma)\setminus\{x_0,z_1,z_2\}$ is not a weakly resolving set of $\Gamma$, the vertex $y_1$ can not resolve $z_1$ and $z_2$, and so $\tilde\partial(y_1,z_2)=\tilde\partial(y_1,z_1)\neq (1,2)$. On the other hand, By Lemma~\ref{arcx1x2}, one has $\tilde\partial(y_1,z_2)=\tilde\partial(z_1,y_2)=(1,2)$. This contradiction shows that $\tilde\partial(y_2,z_2)=(1,1)$. Noting that $X_2=\{y_2,z_2\}$, we get the desired result.
$\qed$

\begin{lemma}\label{arc11}
  If $|X_2|\geq 2$, then $\tilde\partial(w_1,w_2)=(1,1)$ for each $w_1\in X_1$ and each $w_2\in X_2$.
\end{lemma}
\proof We claim that $\partial(x_1,w_2)=1$. In fact, if $|X_1|=1$, then $x_1$ is the unique vertex at distance $1$ from $x_0$, which implies that $\partial(x_1,w_2)=1$ since $\partial(x_0,w_2)=2$. If $|X_2|\geq 2$, then by Lemma~\ref{arcx1x2}, one gets $\tilde\partial(x_1,w_2)=\tilde\partial(x_1,x_2)=(1,\partial(x_2,x_1))$, and so $\partial(x_1,w_2)=1$. Hence, the claim is valid.

Choose a vertex $y_2\in X_2\setminus\{w_2\}$. It follows from Lemma~\ref{x2clique} that
$\tilde\partial(w_2,y_2)=(1,1)$. By the claim, one has $\tilde\partial(w_2,x_1)=(\partial(w_2,x_1),1)$. It follows from Lemma \ref{x1x2} (ii) that $\tilde\partial(w_2,x_0)=(1,2)$. By Lemma~\ref{3vertices}~(i), we get $\tilde\partial(w_2,x_1)=\tilde\partial(w_2,y_2)=(1,1)$.
If $|X_1|=1$, then $w_1=x_1$, and so $\tilde\partial(w_1,w_2)=(1,1)$. If $|X_1|\geq 2$, the desired result follows from Lemma~\ref{arcx1x2}.
$\qed$

\begin{prop}\label{case4}
  If $|X_2|\geq 2$, then $\Gamma$ is isomorphic to
  $$
  G_1[K_1,K_{t},K_{n-t-1}]\qquad \text{or} \qquad  G_2[K_1,K_{t},K_{n-t-1}],
  $$
  where $1\leq t\leq n-3$.
\end{prop}
\proof By Lemmas~\ref{x1x2},~\ref{x2clique} and~\ref{arc11}, we only need to show that $X_1$ is a clique.
Suppose for the contrary that there exist two distinct vertices $y_1,z_1\in X_1$ such that $\tilde\partial(y_1,z_1)\neq (1,1)$. By Lemma~\ref{arc11}, the vertex $y_1$ resolves $z_1$ and $x_2$.
Any vertex $w$ in $X_2$ with $w\neq x_2$ resolves $\{x_0,x_2\}$ and $\{x_0,z_1\}$, so $V(\Gamma)\setminus\{x_0,z_1,x_2\}$ is a weakly resolving set, a contradiction.
$\qed$

Combining Lemma~\ref{check} and Propositions~\ref{case1},~\ref{case2}~\ref{case3},~\ref{case4}, we get Theorem~\ref{n-2main}.

\section*{Acknowledgements}
Feng is supported by the National Natural
Science Foundation of China (11701281), the Natural Science Foundation of Jiangsu Province (BK20170817)
and the Grant of China Postdoctoral Science Foundation. Wang is supported by the National Natural
Science Foundation of China (12071039). Yang is supported by the Fundamental Research Funds for the Central Universities (2652019319).

\end{CJK*}

\end{document}